\documentclass[12pt]{amsart}
\usepackage{amssymb}
\usepackage{url}
\usepackage{mathtools}
\usepackage{color}
\usepackage[pdftex,lmargin=1in,rmargin=1in,tmargin=1.25in,bmargin=1.25in]{geometry}
\usepackage[bookmarks=true, bookmarksopen=true,%
bookmarksdepth=3,bookmarksopenlevel=2,%
colorlinks=true,%
linkcolor=blue,%
citecolor=blue,%
filecolor=blue,%
menucolor=blue,%
pagecolor=blue,%
urlcolor=blue]{hyperref}

\usepackage[percent]{overpic}

\newread\testin

\makeatletter
\def\input@path{{}{draws/}}
\makeatother

\makeatletter
\newcommand\mi@kern[1]{%
  \settowidth\@tempdima{$\mi@obj^{#1}$}
  \kern-\@tempdima
  #1
  \settowidth\@tempdima{$\mi@obj$}
  \kern\@tempdima
}

\newtoks\mi@toksp
\newtoks\mi@toksb
\newcommand{\manyindices}[5]{
  \def\mi@obj{#5}
  \mi@toksp\expandafter{\mi@kern{#2}}
  \mi@toksb\expandafter{\mi@kern{#1}}
  \@mathmeasure4\textstyle{#5_{#1}^{#2}}
  \@mathmeasure6\textstyle{#5_{#3}^{#4}}
  \dimen0-\wd6 \advance\dimen0\wd4
  \@mathmeasure8\textstyle{\hphantom{{}_{#1}^{#2}}#5^{\the\mi@toksp#4}_{\the\mi@toksb#3}}
  \hbox to \dimen0{}{\kern-\dimen0\box8}
}
\makeatother 

\newcommand{\RR}{\mathbb R}

\newcommand{\bD}{\mathbb{D}}

\newcommand{\co}{\colon}

\newcommand{\bdy}{\partial}

\newcommand{\lbracket}{[}
\newcommand{\rbracket}{]}

\DeclareMathOperator{\Sym}{Sym}

\DeclareMathOperator{\ind}{ind}

\DeclareMathOperator{\br}{br} 

\theoremstyle{plain}

\theoremstyle{definition}

\theoremstyle{remark}
\newtheorem{example}{Example}

\newcommand\HH{\mathit{HH}}
\newcommand\Hochschild\HH

\newcommand{\alphas}{{\boldsymbol{\alpha}}}
\newcommand{\betas}{{\boldsymbol{\beta}}}

\makeatletter
\newcommand\honestalg[3]{\bigl\lbracket
\begin{smallmatrix} #1\@ifempty{#3}{}{&#3} \\ #2 \end{smallmatrix}
\bigr\rbracket}

\makeatother

\newcommand{\cD}{\mathcal{D}}

\newcommand{\vx}{{\vec{x}}}
\newcommand{\vy}{{\vec{y}}}

\begin{document}
\title{Errata to ``A cylindrical reformulation of Heegaard Floer homology''}

\author{Robert Lipshitz}
\address{Department of Mathematics, Columbia University\\
  New York, NY 10027}
\date{\today}
\email{lipshitz@math.columbia.edu}
\maketitle

This note has two parts. The first part explains a serious gap in the
proof of the index formula in~\cite[Section
4]{Lipshitz06:CylindricalHF}, discovered by John Pardon. We explain
the gap in Section~\ref{sec:ind-gap} and how to correct the proof of
the index formula in
Section~\ref{sec:ind-fix}. Section~\ref{sec:smaller} acknowledges and
corrects four smaller errors, not affecting the main results
of~\cite{Lipshitz06:CylindricalHF}.

\subsection*{Acknowledgments} I thank John Pardon, Clifford Taubes
and Guangbo Xu for pointing out errors
in~\cite{Lipshitz06:CylindricalHF}, and for helpful conversations
about how to correct these errors. I also thank John Pardon for
helpful comments on drafts of these errata.

\section{The index formula for embedded curves}
\subsection{The gap}\label{sec:ind-gap}
\subsubsection{What is correct}
In the cylindrical formulation, there are two steps to studying the
expected dimensions of the moduli spaces. The first step is to
consider the $\overline{\bdy}$-operator for maps
\begin{equation}\label{eq:map}
u\co (S,\bdy S)\to \bigl(\Sigma\times[0,1]\times\RR,
(\alphas\times\{1\}\times\RR)\cup (\betas\times\{0\}\times\RR)\bigr)
\end{equation}
for a fixed homeomorphism type of $S$. It is shown that the index of
the $\overline{\bdy}$-operator for such maps is given by
\[
\ind(u)=g-\chi(S)+2e(A)
\]
where $g$ is the genus of the Heegaard surface (or, more importantly,
the number of negative (equivalently positive) ends of $u$), and
$e(A)$ is the Euler measure of the domain $A$ in $\Sigma$ of the map
$u$. (This formula holds whether or not $u$ is holomorphic.)

The cylindrical formulation of Heegaard Floer homology corresponds to counting embedded holomorphic
curves of the form~\eqref{eq:map}. So, the second step in studying the
index is to show that for embedded curves, $\chi(S)$ is determined by
the homology class $A$. It is shown in~\cite[Proposition~4.2 and
Corollary~4.3]{Lipshitz06:CylindricalHF} that at an embedded
holomorphic curve, $\chi(S)$ is given by
\begin{align}
\chi(S)&=g-n_\vx(A)-n_\vy(A)+e(A)\label{eq:chi}\\
\shortintertext{so}
\ind(u)&=e(A)+n_\vx(A)+n_\vy(A).\label{eq:ind}
\end{align}

The proofs in~\cite{Lipshitz06:CylindricalHF} of
Formulas~\eqref{eq:chi} and~\eqref{eq:ind} at an embedded holomorphic
curve $u$, with respect to any almost complex structure satisfying the
conditions~\cite[(J1)--(J5), p.\ 959]{Lipshitz06:CylindricalHF})
(including non-generic almost complex structures of this form), are
correct.

Homology classes of curves in $\Sigma\times[0,1]\times\RR$ correspond
to homotopy classes of disks in the symmetric product.  If $A$ is
represented by an embedded holomorphic curve with respect to the
product complex structure on $\Sigma\times[0,1]\times\RR$, it follows
from the tautological correspondence that $\ind(u)$ agrees with the
Maslov index in the symmetric product. So, in these cases,
Formula~\eqref{eq:ind} computes the Maslov index for disks in
$\Sym^g(\Sigma)$.

\subsubsection{What more one wants}

It is natural to be interested in the index at homology classes not
represented by embedded holomorphic curves, for two reasons:
\begin{enumerate}
\item\label{item:ind-adds} One wants to know that the right-hand side of
  Formula~\eqref{eq:ind} is additive, so one can use it to define or
  compute the relative grading on the Heegaard Floer complexes.
\item It is tidier to know that Formula~\eqref{eq:ind} always agrees
  with the Maslov index in $\Sym^g(\Sigma)$; the Maslov index is
  defined whether or not there is a holomorphic representative.
\end{enumerate}
 
Note that S.~Sarkar has given a combinatorial proof that
Formula~\eqref{eq:ind} is additive, in the process of generalizing it
to give a formula for the Maslov index of higher holomorphic
polygons~\cite{Sarkar06:IndexTriangles}.

\subsubsection{What is wrong}
To generalize Formula~\eqref{eq:chi} to homology classes not admitting
holomorphic representatives, we need some class of maps $u$ which is
broader than holomorphic maps but for which $\chi(S)$ is still
determined. To show that the right-hand side of Formula~\eqref{eq:ind}
agrees with the Maslov index in $\Sym^g(\Sigma)$ for homology classes
without holomorphic representatives, we also want these maps $u$ to
correspond to disks in $\Sym^g(\Sigma)$. Such classes of maps $u$ were
proposed in~\cite[Lemmas 4.1 and 4.9]{Lipshitz06:CylindricalHF}:

\noindent\textbf{Lemma 4.1}\ \emph{
  Suppose \(A\in\pi_2(\vx,\vy)\) is a positive homology class.
  Then there is a Riemann surface
  with boundary and corners \(\overline{S}\) and smooth map \(u\co S\to
  \Sigma\times[0,1]\times\RR\) (where \(S\) denotes the complement in \(\overline{S}\) of the
  corners of \(\overline{S}\)) in the homology class \(A\) such that 
  \begin{enumerate}
  \item \(u^{-1}(C_\alpha\cup C_\beta)=\bdy S\).
  \item  For each \(i\), \(u^{-1}(\alpha_i\times\{1\}\times\RR)\) and
    \(u^{-1}(\beta_i\times\{0\}\times\RR)\) each consists of one arc in
    \(\bdy S\).
  \item The map \(u\) is \(J\)--holomorphic in a
    neighborhood of \((\pi_\Sigma\circ u)^{-1}(\alphas\cup\betas)\) for some
    \(J\) satisfying (\textbf{J1})--(\textbf{J5}) (in fact, for
    \(j_\Sigma\times j_\bD\)).
  \item For each component of \(S\), either
    \begin{itemize}
    \item The component is a disk with two boundary punctures and the map is
      a diffeomorphism to \(\{x_i\}\times[0,1]\times\RR\) for some
      \(x_i\in\alphas\cap\betas\) (such a component is a \emph{degenerate disk}) or
    \item The map \(\pi_\Sigma\circ u\) extends to a branched covering map
      \(\overline{\pi_\Sigma\circ u}\), none
      of whose branch points map to points in \(\alphas\cap\betas\).
    \end{itemize}
  \item All the corners of \(S\) are acute.
  \item \label{item:last}The map \(u\) is an embedding.
  \end{enumerate}
}

\noindent\textbf{Lemma 4.9}\ \emph{
  Suppose \(A\) is a positive homology class.  Then we can
  represent \(A+[\Sigma]\) by a map \(u\co S\to \Sigma\times[0,1]\times\RR\) satisfying all the
  conditions of Lemma~4.1 and such that,
  additionally,
  \begin{itemize}
  \item The map \(\pi_\bD\circ u\) is a \(g\)--fold branched covering
    map with
    all its branch points of order \(2\)
  \item The map \(u\) is holomorphic near the preimages of the branch
    points of \(\pi_\bD\circ u\).
  \end{itemize}  
}
  
The proof of Lemma 4.1 has two gaps:
\begin{enumerate}
\item In the proof, one starts by gluing up the domain of $u$ to
  produce a surface. One wants to ensure that the only corners
  correspond to the points in $\vx\cup \vy$. The proof says to start
  with any maximal gluing and then make some local changes, but is
  imprecise or incorrect about how to do so.
\item The argument for ensuring that the map is an embedding
  (property~(\ref{item:last})) is incorrect. First, some map, not
  necessarily an embedding is constructed. Then, the proof says:
  ``Modifying \(S_1\) and \(p_{\Sigma,1}\times p_{\bD,1}\) near the
  double points of \(p_{\Sigma,1}\times p_{\bD,1}\) we can obtain a
  new map \(u\co S\to \Sigma\times[0,1]\times\RR\) satisfying all of the stated
  properties''. Typically, such a modification is not possible while
  keeping $\pi_\Sigma\circ u$ a branched map; see
  Example~\ref{eg:cant-unbranch}.
\end{enumerate}
The proof of Lemma 4.9 builds on Lemma 4.1, and has the same gaps.

As we will see, the first point can be resolved by being more careful
in the construction, following~\cite[Lemma
2.17]{OS04:HolomorphicDisks} (see also~\cite[Lemma 10.3]{LOT1}).  The
second point is more serious, as the following example (explained to
me by J.~Pardon) shows.

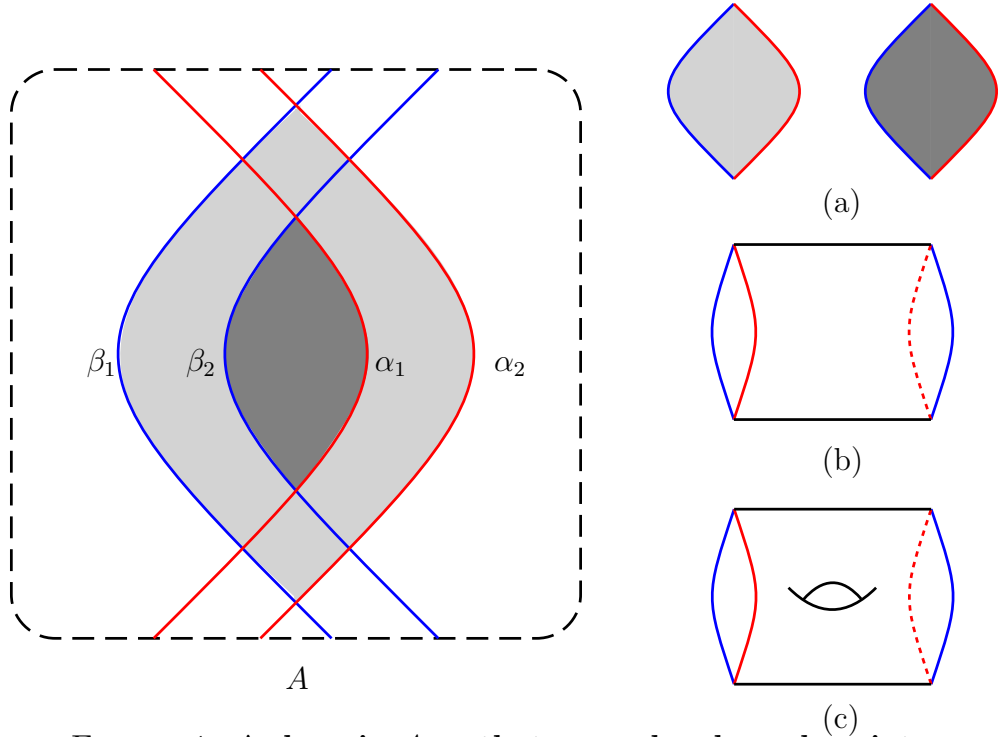
\begin{figure}
  \centering
  \begin{overpic}[tics=10,width=.8\textwidth]{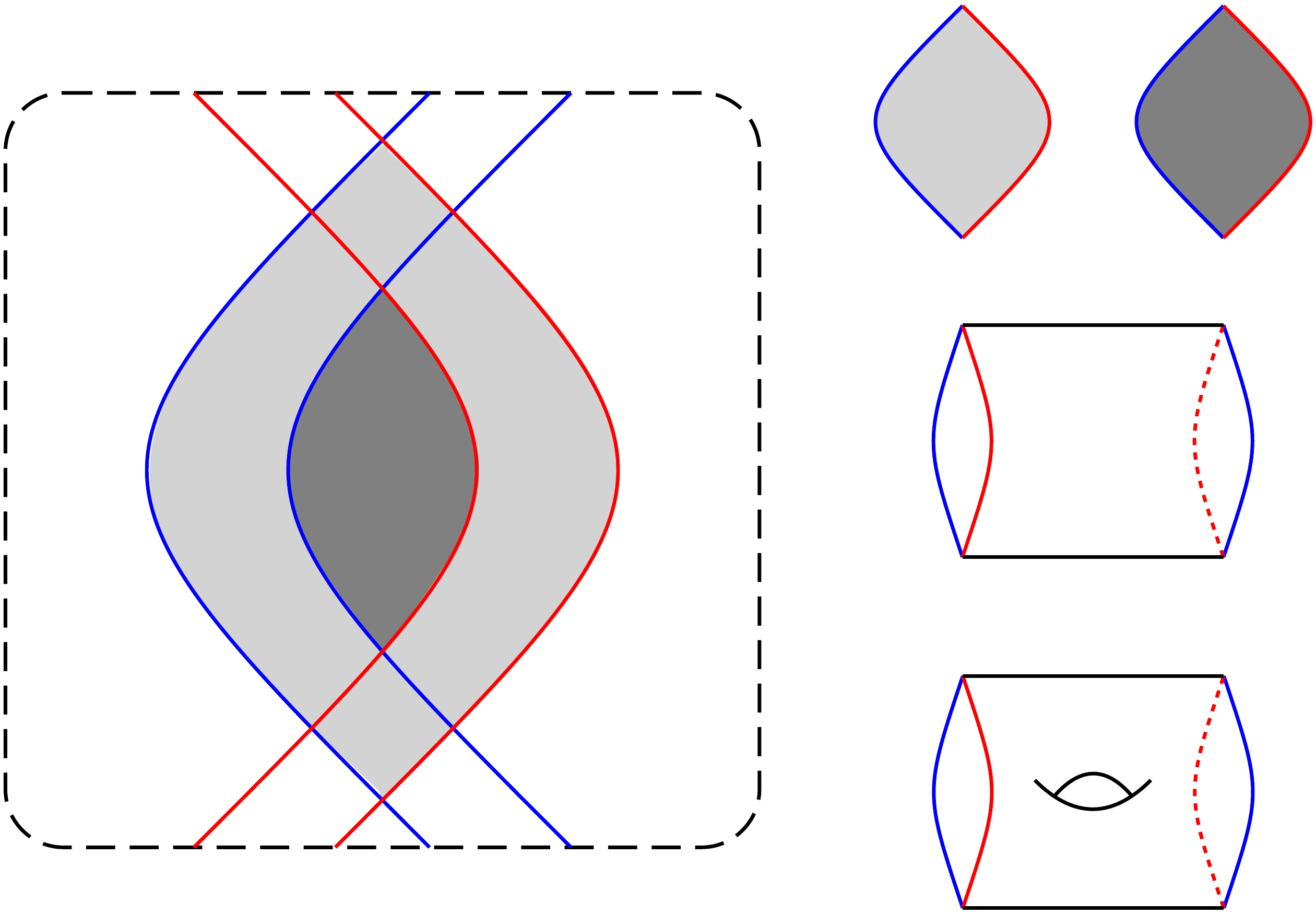}
    \put(28,0){$A$}
    \put(8,32){$\beta_1$}
    \put(18,32){$\beta_2$}
    \put(37,32){$\alpha_1$}
    \put(49,32){$\alpha_2$}
    \put(82,48){(a)}
    \put(82,22){(b)}
    \put(82,-4){(c)}
  \end{overpic}
  \caption{\textbf{A domain $A$ so that $\pi_\Sigma\circ u$ has branch points.}}
  \label{fig:br-domain}
\end{figure}

\begin{example}\label{eg:cant-unbranch}
  Consider the domain $A$ shown in Figure~\ref{fig:br-domain}. There
  is an obvious holomorphic representative $S\to
  \Sigma\times[0,1]\times\RR$ where $S$ is the disjoint union of two
  disks (bigons). This representative has a positive
  double-point. Resolving the double point gives a map
  $S'\to\Sigma\times[0,1]\times\RR$ where $S'$ is an annulus. Indeed,
  Formula~\eqref{eq:chi} predicts the embedded Euler
  characteristic $\chi=2-6/4-6/4+1=0$.

  But we can also find other, non-holomorphic representatives of this
  domain. For example, take $S$ to be a surface of genus $1$ with $2$ boundary
  components. Then we can find branched maps $S\to\Sigma$ and $S\to
  \bD$ representing the domain $A$. It is easy to arrange this map to
  satisfy the conditions in Lemma 4.1 except for being an
  embedding. Resolving double points decreases the Euler
  characteristic of $S$, which is already lower than the Euler
  characteristic predicted by Formula~\ref{eq:chi}; so, if we could
  resolve them (without losing the other properties in Lemma~4.1),
  this would contradict Proposition~4.2.
\end{example}

\subsection{Revised proofs of the main results}\label{sec:ind-fix}
We can salvage the main result by weakening the conditions in Lemmas
4.1 and 4.9 to allow $u$ to have double points, and strengthening
Proposition 4.2, Corollary 4.3 and the proofs of Propositions 4.8 and
Corollary 4.10 to curves with double points.

The revised Lemma 4.1 reads as follows:

\noindent\textbf{Lemma 4.1$\boldsymbol{'}$}\ \emph{
  Suppose \(A\in\pi_2(\vx,\vy)\) is a positive homology class.
  Then there is a Riemann surface
  with boundary and corners \(\overline{S}\) and smooth map \(u\co S\to
  \Sigma\times[0,1]\times\RR\) (where \(S\) denotes the complement in \(\overline{S}\) of the 
  corners of \(\overline{S}\)) in the homology class \(A\) such that 
  \begin{enumerate}
  \item \label{item:4.1p:bdy}\(u^{-1}(C_\alpha\cup C_\beta)=\bdy S\).
  \item \label{item:4.1p:one-arc} For each \(i\), \(u^{-1}(\alpha_i\times\{1\}\times\RR)\) and
    \(u^{-1}(\beta_i\times\{0\}\times\RR)\) each consists of one arc in
    \(\bdy S\).
  \item\label{item:4.1p:holo} The map \(u\) is \(J\)--holomorphic in a
    neighborhood of \(\bdy S\)
    for some
    \(J\) satisfying (\textbf{J1})--(\textbf{J5}) (in fact, for
    \(j_\Sigma\times j_\bD\)).%
  \item\label{item:4.1p:components} For each component of \(S\), either
    \begin{itemize}
    \item The component is a disk with two boundary punctures and the map is
      a diffeomorphism to \(\{x_i\}\times[0,1]\times\RR\) for some
      \(x_i\in\alphas\cap\betas\) (such a component is a \emph{degenerate disk}) or
    \item The map \(\pi_\Sigma\circ u\) extends to a branched covering map
      \(\overline{\pi_\Sigma\circ u}\), none
      of whose branch points map to points in \(\alphas\cap\betas\).
    \end{itemize}
  \item \label{item:4.1p:acute} All the corners of \(S\) are acute.
  \item \label{item:4.1p:double-pts}The map \(u\) has at worst transverse double point singularities.
  \end{enumerate}
}

(For convenience, we have also weakened
Condition~(\ref{item:4.1p:holo}); the resulting condition is sufficient
for the other results to go through and requires one fewer step to
achieve.)

The statement of Lemma 4.9 does not need any revisions, except that
``all the conditions of Lemma~4.1'' now refers to Lemma $4.1'$; and we
should have assumed that $g>1$:

\noindent\textbf{Lemma 4.9$\boldsymbol{'}$}\ \emph{Assume that $g>1$.
  Suppose \(A\) is a positive homology class.  Then we can
  represent \(A+[\Sigma]\) by a map \(u\co S\to \Sigma\times[0,1]\times\RR\) satisfying all the
  conditions of Lemma~4.1$'$ and such that,
  additionally,
  \begin{itemize}
  \item The map \(\pi_\bD\circ u\) is a \(g\)--fold branched covering
    map with
    all its branch points of order \(2\)
  \item The map \(u\) is holomorphic near the preimages of the branch
  points of \(\pi_\bD\circ u\).
\end{itemize}
}

Proposition 4.2 now reads:

\noindent\textbf{Proposition 4.2$\boldsymbol{'}$}\ \emph{
  Let \(u\co S\to \Sigma\times[0,1]\times\RR\) be a map satisfying the conditions
  enumerated in the Lemma 4.1',
  representing a homology class \(A\).  Suppose that $u$ has $d_+$
  positive double points and $d_-$ negative double points. Then the Euler characteristic
  \(\chi(S)\) is given by  
  \[
    \chi(S)=g-n_\vx(A)-n_\vy(A)+e(A)+2(d_+-d_-).
  \]
}

With these changes, it is clearer to state Proposition 4.8 as follows:

\noindent\textbf{Proposition 4.8$\boldsymbol{'}$}\ \emph{
The Maslov index (in the symmetric product) of a domain $A\in\pi_2(\vx,\vy)$ is given by
\begin{equation}\label{eq:mu-equals}
\mu(A)=e(A)+n_\vx(A)+n_\vy(A).
\end{equation}
This agrees with the index in the cylindrical setting at any embedded
holomorphic curve (with respect to an almost complex structure
satisfying conditions (J1)--(J5)).
}

\begin{proof}[Proof of Lemma $4.1'$]
  This construction is adapted from the proof of~\cite[Lemma
    2.17]{OS04:HolomorphicDisks}. Let $\{\cD_i\}$ denote the
  components of $\Sigma\setminus(\alphas\cup\betas)$. Write $A=\sum
  n_i\cD_i$; assume that we have ordered the $\cD_i$ so that if $i<j$
  then $n_i\leq n_j$. Build a surface $S_0$ by taking, for each $i$,
  $n_i$ copies of $\cD_i$; denote these copies $\cD_i^{(j)}$. Glue
  these together as follows:
  \begin{itemize}
  \item If $\cD_i$ and $\cD_j$ ($i<j$) share a common $\alpha$-arc $a$
    then for each $k=1,\dots,n_i$ glue $\cD_i^{(k)}$ to $\cD_j^{(k+n_j-n_i)}$ along $a$.
  \item If $\cD_i$ and $\cD_j$ ($i<j$) share a common $\beta$-arc $b$
    then for each $k=1,\dots,n_i$ glue $\cD_i^{(k)}$ to $\cD_j^{(k)}$ along $b$.
  \end{itemize}
  The resulting surface $S_0$ comes equipped with a map
  $u_{\Sigma,0}\co S_0\to\Sigma$. The surface $S_0$ is obviously a
  smooth surface-with-boundary away from
  $u_{\Sigma,0}^{-1}(\alphas\cap\betas)$. Next, consider the behavior
  of $S_0$ near a point $p\in\alphas\cap\betas$. Let
  $\cD_{i_1},\dots,\cD_{i_4}$ be the four regions incident to
  $p$. (Some of the $D_{i_j}$ might be the same.) If
  $p\notin\vx\cup\vy$ or $p\in\vx\cap\vy$ then the coefficients of
  the $D_{i_j}$ in $A$ have the form $n,n+k,n+\ell,n+k+\ell$ for some
  $n,k,\ell\geq 0$. Reordering the $i_j$ if necessary, assume that $D_{i_1}$ and $D_{i_2}$ share a common $\alpha$-arc. Then the following $D_{i_j}^{(m)}$ are glued together:
  \begin{align*}
  \{D_{i_1}^{(1)},D_{i_2}^{(1+k)},D_{i_3}^{(1)},D_{i_4}^{(1+k)}\},\dots,
  \{D_{i_1}^{(n)},D_{i_2}^{(n+k)},D_{i_3}^{(n)},D_{i_4}^{(n+k)}\}&,\\
  \{D_{i_2}^{(1)},D_{i_4}^{(1)}\},\dots,\{D_{i_2}^{(k)},D_{i_4}^{(k)}\}&\\
  \{D_{i_3}^{(n+1)},D_{i_4}^{(n+k+1)}\},\dots,\{D_{i_3}^{(n+\ell)},D_{i_4}^{(n+k+\ell)}\}&.
  \end{align*}
  In particular, near each point in $u_{\Sigma,0}^{-1}(p)$, $S_0$ is
  again a smooth surface-with-boundary and the map $u_{\Sigma,0}$ is a
  homeomorphism onto its image.

  Now, suppose that $p\in\vx\setminus\vy$ or
  $p\in\vy\setminus\vx$. Then the coefficients of $A$ near $p$ can be
  written as one of $\{a+1,a+k,a+\ell,a+k+\ell\}$,
  $\{a,a+k+1,a+\ell,a+k+\ell\}$, $\{a,a+k,a+\ell+1,a+k+\ell\}$ or
  $\{a,a+k,a+\ell,a+k+\ell+1\}$. In the first case, the glued regions are
  \begin{align*}
    \{D_{i_1}^{(a+1)},D_{i_2}^{(a+k)},D_{i_4}^{(a+k)},D_{i_3}^{(a)},D_{i_1}^{(a)},D_{i_2}^{(a+k-1)},D_{i_4}^{(a+k-1)},D_{i_3}^{(a-1)},\dots,D_{i_1}^{(1)}\}&\\
    \{D_{i_2}^{(1)},D_{i_4}^{(1)}\},\dots,\{D_{i_2}^{(k-1)},D_{i_4}^{(1)}\}&\\
    \{D_{i_3}^{(a+1)},D_{i_4}^{(a+k+1)}\},\dots,\{D_{i_3}^{(a+\ell)},D_{i_4}^{(a+k+\ell})\}&.
  \end{align*}
  In particular, the preimage of $p$ consists of $(k-1)+\ell$
  preimages which are smooth boundary points, and near which
  $u_{\Sigma,0}$ is a homeomorphism onto its image; and one preimage
  which looks like a boundary branch point. Call this last preimage a
  \emph{bad point}. If we choose a smooth structure on $S_0$ making
  the bad point a $\pi/2$ corner then the map $u_{\Sigma,0}$ is of the
  form $z\mapsto z^{4a+1}$ near this point.

  The other three cases are similar, in that all but one of the
  preimages of $p$ lie on the smooth boundary of $S$, and near them
  the map $u_{\Sigma,0}$ is a local homeomorphism; and there is one
  remaining, \emph{bad point} near which we can straighten $S_0$ and
  view $u_{\Sigma,0}$ as a branched map.
  In particular, for each $\alpha_i$, $u_{\Sigma,0}^{-1}(\alpha_i)\cap \bdy
  S_0$ consists of a union of some circles $C^\alpha_{i,j}$ and
  possibly a single arc $A^\alpha_i$; and similarly for each
  $\beta_i$.

  The surface $S_0$ has corners, which are in bijective correspondence
  with
  $(\vx\cup\vy)\setminus(\vx\cap \vy)$. The map $u_{\Sigma,0}$ may have
  branch points at some of these corners, say $p_1,\dots,p_k$. If
  $p_i$ has total angle $n\pi/2$, make $(n-1)/2$ cuts in $S_0$ at $p_i$,
  as in Figure~\ref{fig:corner-cuts}. Let $S_1$ be the resulting
  surface and $u_{\Sigma,1}$ the resulting map to $\Sigma$.

  \begin{figure}
    \centering
    \begin{overpic}[tics=10,scale=1]{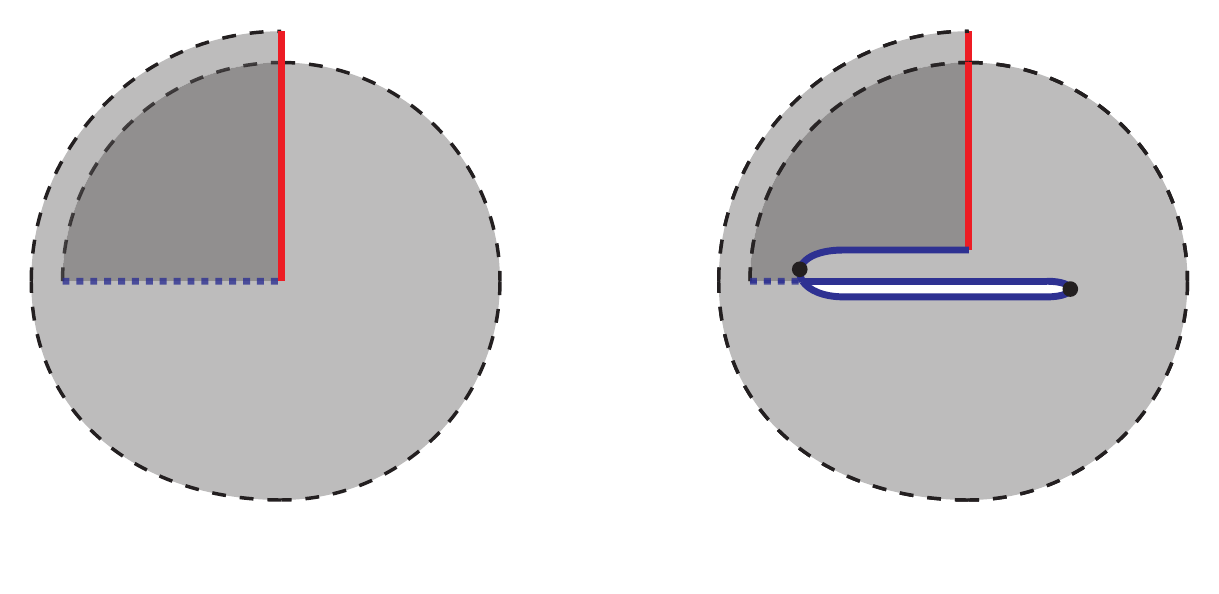}
      \put(21,0){$S_0$}
      \put(78,0){$S_1$}
      \put(24,26){$x_i$}
      \put(80.5,28){$x_i$}
    \end{overpic}
    \caption{\textbf{Making cuts at the corners.} The figure shows a
      region of $S_0$ (left) and $S_1$ (right); in $S_0$ there is a
      branch point at $x_i$. The darker region is covered with
      multiplicity $2$. On the right, the two dots are boundary branch
      points. We made cuts along the $\beta$-arcs; we could equally
      well have made cuts along the
      $\alpha$-arcs.\label{fig:corner-cuts}}
  \end{figure}
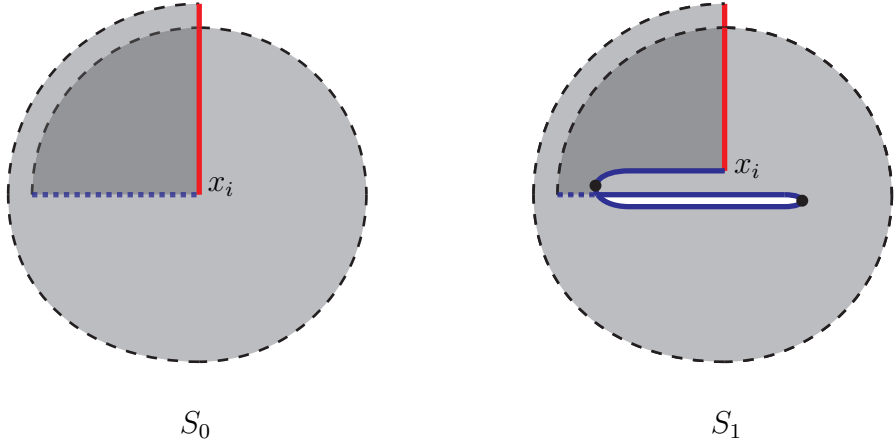

  Next, we modify $(u_{\Sigma,1},S_1)$ to a new surface whose corners
  correspond to $\vx\amalg \vy$; that is, we introduce corners corresponding to
  points in $\vx\cap \vy$. For each point $x_i\in \vx\cap\vy$, if $x_i$ is
  disjoint from $u_{\Sigma,1}(\bdy S_1)$ then take the disjoint union
  of $S_1$ with a twice punctured disk, and define $u_{\Sigma_2}$ to
  map the twice punctured disk by a constant map to $x_i$. If $x_i$ is
  not disjoint from $u_{\Sigma,1}(\bdy S_1)$ then choose an arc in
  $\bdy S_1$ covering $x_i$ and make a small slit in the arc starting
  at $x_i$. (This introduces two new corners, both mapping to $x_i$,
  and a boundary branch point.) See Figure~\ref{fig:degen-corners}.
  After this modification, the surface has exactly $2g$
  corners, corresponding to $\vx\amalg\vy$. Call the result
  $(S_2,u_{\Sigma,2})$.

  \begin{figure}
    \centering
    \begin{overpic}[tics=20,height=1.5in]{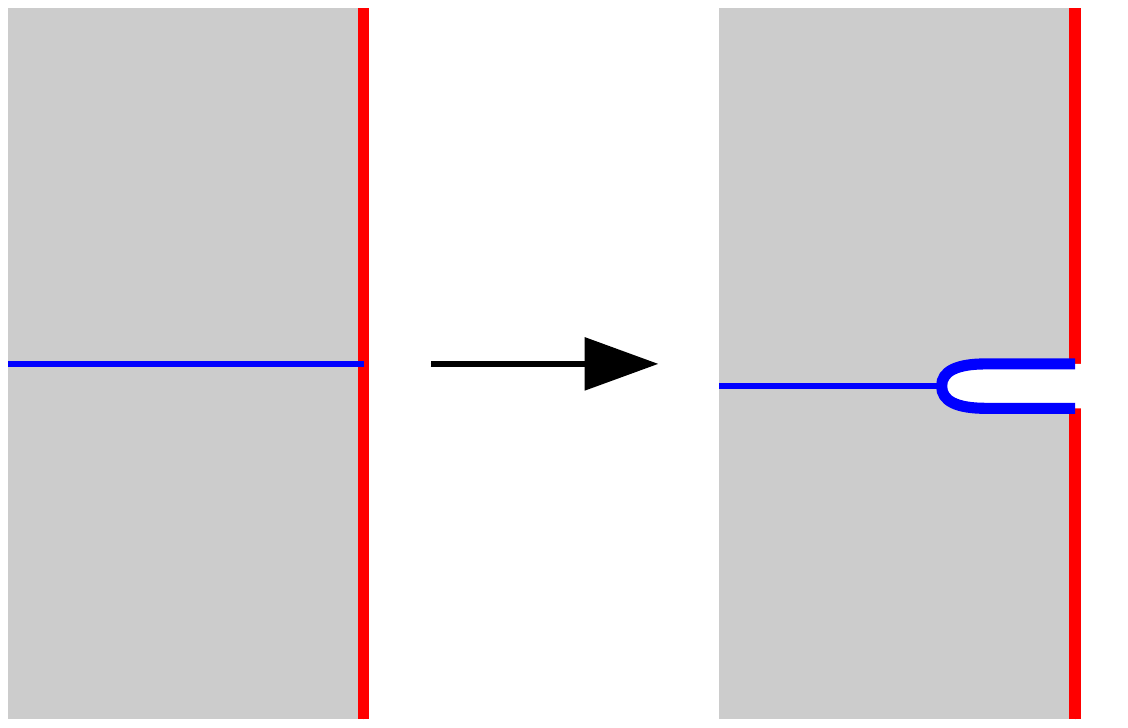}
      \put(34,27){$x_i$}
    \end{overpic}
    \caption{\textbf{Adding slits at degenerate corners.}}
    \label{fig:degen-corners}
  \end{figure}

  \begin{figure}
    \centering
    \begin{overpic}[tics=20,width=\textwidth]{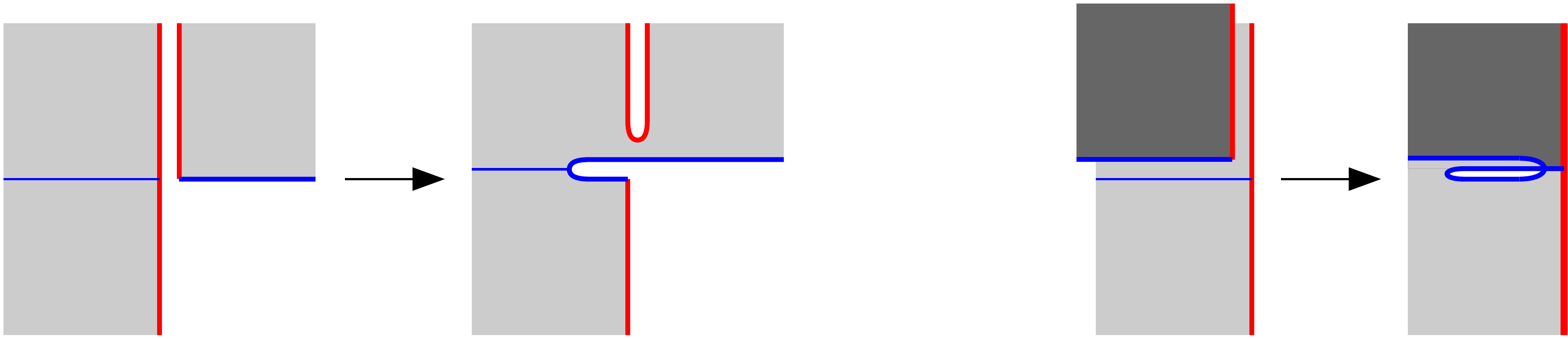}
    \end{overpic}
    \caption{\textbf{Splicing a corner and an edge.}}
    \label{fig:splice-corners}
  \end{figure}

  Next, we modify $(S_2,u_{\Sigma,2})$ to a new pair
  $(S_3,u_{\Sigma,3})$ so that for each $i$,
  $u_{\Sigma,3}^{-1}(\alpha_i)\cap\bdy S_3$ consists of a single arc
  (and no circles); and similarly for each $\beta$-circle.  In the
  process, we will introduce some more boundary branch points.
  Suppose that $C$ is a boundary component of $S_2$ which is mapped
  entirely to $\alpha_i$. Let $x_i\in \vx$ be the corner on
  $\alpha_i\cap \beta_{j}$ (for some $j$) and $p_i$ the corresponding
  corner of $S_2$. Make a small slit in $S_2$ along
  $u_{\Sigma,2}^{-1}(\beta_j)$ starting at $C\cap
  u_{\Sigma,2}^{-1}(x_i)$, and glue one edge of the resulting surface
  to the $\alpha$- or $\beta$-arc near $p_j$ in such a way that $u_{\Sigma,2}$
  induces a branched map from the result. (There are two cases for the
  local geometry here; see Figure~\ref{fig:splice-corners}.) This
  reduces the number of boundary components of $S$ mapped to
  $\alpha_i$ by $1$; repeat for the other $\alpha$-boundary circles of
  $S_2$.  Modify boundary components mapped entirely to $\beta_i$
  similarly. Call the result $(S_3,u_{\Sigma,3})$; this pair has the
  property that $u_{\Sigma,3}^{-1}(\alpha_i)\cap\bdy S_3$ consists of
  a single arc (and no circles).

  The map $u_{\Sigma,3}$ and the complex structure on $\Sigma$ induce
  a complex structure on $S_3$. Let $U$ denote a tubular neighborhood
  of $\bdy S_3$. Choose a holomorphic map $u_{\bD,3}\co U\to
  [0,1]\times\RR$ so that:
  \begin{itemize}
  \item $u_{\bD,3}$ sends each $\alpha$-arc in $\bdy S_3$ to
    $\{1\}\times\RR$ and each $\beta$-arc to $\{0\}\times\RR$.
  \item Near each corner of $S_3$ corresponding to a point in $\vx$,
    $u_{\bD,3}$ is asymptotic to $-\infty$.
  \item Near each corner of $S_3$ corresponding to a point in $\vy$,
    $u_{\bD,3}$ is asymptotic to $+\infty$.
  \item The map $u_{\bD,3}$ is a local diffeomorphism (i.e., has
    non-vanishing derivative).
  \end{itemize}

  Extend $u_{\Sigma,3}$ arbitrarily to the rest of $S_3$. Then
  $u_{\Sigma,3}\times u_{\bD,3}$ is a map to
  $\Sigma\times[0,1]\times\RR$. By construction, this map satisfies
  Conditions~(\ref{item:4.1p:bdy}),~(\ref{item:4.1p:one-arc}),~(\ref{item:4.1p:holo})~(\ref{item:4.1p:components})
  and~(\ref{item:4.1p:acute}).
  Perturbing $u_{\Sigma,3}\times u_{\bD,3}$ slightly (without changing
  it near the boundary) gives a map $u\co S\to
  \Sigma\times[0,1]\times\RR$ satisfying
  Condition~(\ref{item:4.1p:double-pts}), as well.
\end{proof}

\begin{proof}[Proof of Proposition $4.2'$]
  The proof is essentially the same as the original proof of
  Proposition~4.2, noting that each double point leads to two
  intersections of $u$ and $u'$. We spell this out.

  First, note that each degenerate disk adds $1$ to $\chi(S)$, $1$ to
  $g$, $0$ to $e(A)$ and $0$ to
  $2(d_+-d_-)-n_{\vx}(A)-n_{\vy}(A)$. Thus, each such disk changes the
  two sides of the formula in identical ways, and so we may assume
  there are no degenerate disks.

  Next, by the Riemann-Hurwitz formula,
  \[
  e(S) = e(A)-\br(\pi_\bD\circ u),
  \]
  where $\br(\pi_\bD\circ u)$ denotes the ramification degree of
  $\pi_\bD\circ u$. (For example, if all branch points of
  $\pi_\bD\circ u$ have order $2$ then $\pi_\bD\circ u$ is just the
  number of branch points.) Moreover, since $S$ has $2g$
  $\pi/2$-corners,
  \[
  \chi(S)=e(S)+g/2;
  \]
  so, we want to compute $\br(\pi_\bD\circ u)$.
    
  Let 
  \[
  \tau_r\co \Sigma\times[0,1]\times\RR\to \Sigma\times[0,1]\times\RR
  \]
  be translation by $r$ units in the $\RR$-direction, i.e.,
  $\tau_r(p,s,t)=(p,s,t+r)$. Let $\partial/\partial t$ denote the
  tangent vector field to $\RR$. Then, for $\epsilon$ sufficiently small, we have
  \begin{align*}
    \br(\pi_\bD\circ u)&=\#\{\text{tangencies of }\partial/\partial
    t\text{ and u}\}\\
    &=\#(u\cap \tau_\epsilon\circ u)-2(d_+-d_-),
  \end{align*}
  where all counts are with multiplicity. (Tangencies along the
  boundary, boundary double points, boundary intersection points and
  boundary branch points each count for $1/2$.) The term $2(d_+-d_-)$
  comes from the fact that each positive (respectively negative)
  double point of $u$ contributes $2$ intersections between $u$ and
  $\tau_\epsilon\circ u$.

  The fact that $u$ is holomorphic near its boundary implies that
  $\#(u\cap \tau_\epsilon\circ u)=\#(u\cap \tau_R(u))$ for any $R\in
  \RR$. When $R$ is sufficiently large, 
  \[
  \#(u\cap \tau_R(u))=n_\vx(A)+n_\vy(A)-g/2.
  \]
  
  Collecting these equalities,
  \begin{align*}
  \chi(S)&=e(S)+g/2\\
  &=e(A)-\br(\pi_\bD\circ u)+g/2\\
  &=e(A) - \#(u\cap \tau_R\circ u)+2(d_+-d_-)+g/2\\
  &=e(A) - n_\vx(A)-n_\vy(A)+g+2(d_+-d_-),
  \end{align*}
  as desired.
\end{proof}

\begin{proof}[Proof of Lemma 4.9$'$]
  Let $(S_3,u_{\Sigma,3})$ be as in the proof of Lemma~4.1$'$ applied
  to the homology class $A$. 
  Build a pair $(S_4,u_{\Sigma,4})$
  representing $A+[\Sigma]$, and with $S_4$ connected, as
  follows. Start with the disjoint union
  $S_3\amalg \Sigma$. 
  Forget each degenerate disk in $S_3$ and instead make three cuts
  starting from the corresponding $x_i\in\Sigma$, two along the
  $\alpha$-circle and one along the $\beta$-circle, as in the left of
  Figure~\ref{fig:cut-Sigma}. 
  At each remaining point $x_i\in\vx$, cut open
  $\Sigma$ in the same way and glue it to the corresponding corner of $S_3$, as shown
  in Figure~\ref{fig:cut-Sigma}. The result is a connected surface
  $S_4$ and map $u_{\Sigma,4}\co S_4\to \Sigma$ representing the
  homology class $A+[\Sigma]$.

  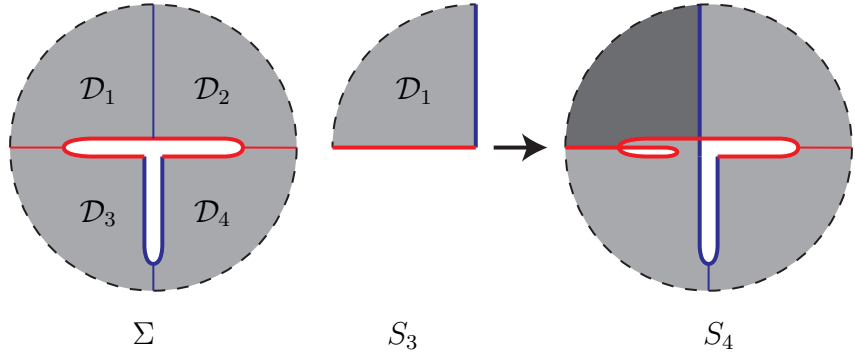
\begin{figure}
    \centering
    \begin{overpic}[tics=10,scale=.75]{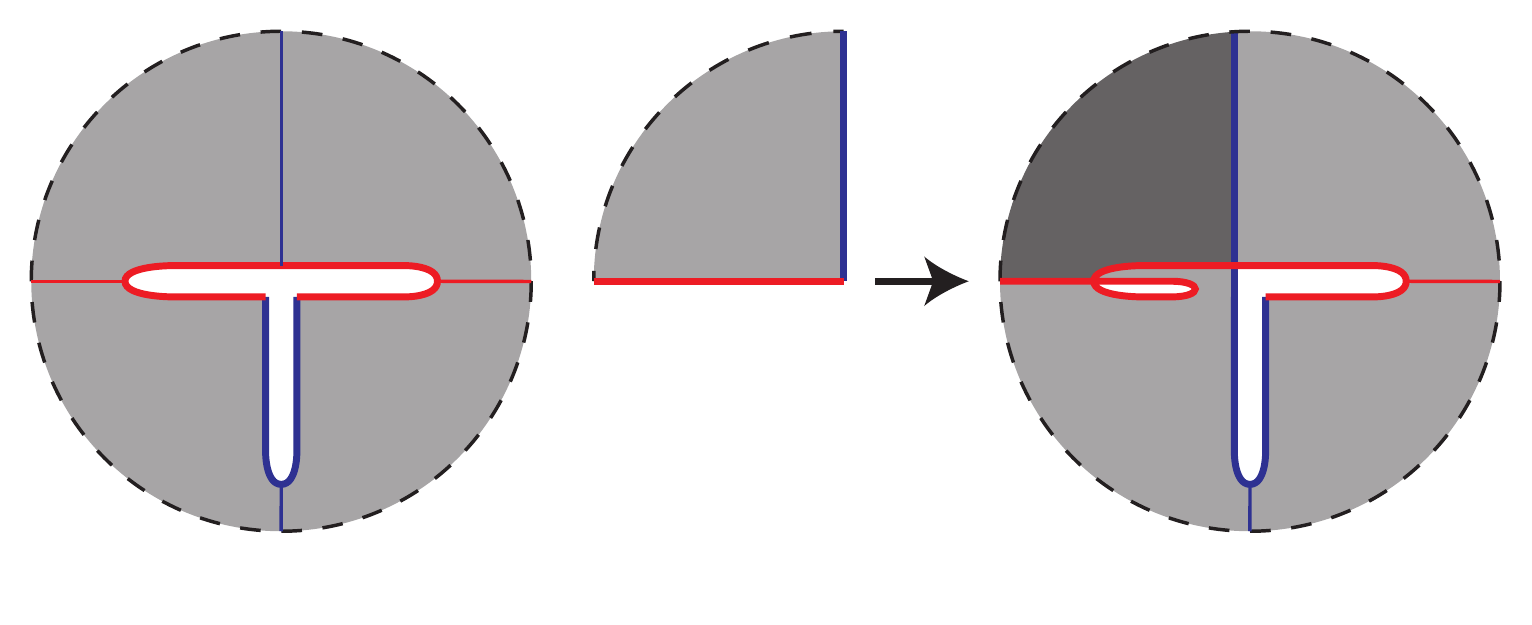}
      \put(16,0){$\Sigma$}
      \put(45,0){$S_3$}
      \put(81,0){$S_4$}
      \put(10,28){$\cD_1$}
      \put(23,28){$\cD_2$}
      \put(10,14){$\cD_3$}
      \put(23,14){$\cD_4$}
      \put(46,28){$\cD_1$}
    \end{overpic}
    \caption{\textbf{Cutting $\Sigma$ and gluing to $S_3$.}}
    \label{fig:cut-Sigma}
  \end{figure}

  We claim that there is a $g$-fold branched covering $u_{\bD,4}\co
  S_4\to [0,1]\times\RR$ sending the $\alpha$-boundary of $S_4$ to
  $\{1\}\times\RR$ and the $\beta$-boundary of $S_4$ to
  $\{0\}\times\RR$. To see this, let $\bdy S_4\times[0,\epsilon)$ be a
  collar neighborhood of $\bdy S_4$ and let $C=\bdy
  S_4\times\{\epsilon/2\}$. Collapsing the circles in $C$ gives a
  surface $S'$ consisting of a union of disks $D$---one for each
  boundary component of $S$---and a closed surface $E$ meeting each
  disk in a single point. If $D_i$ has $2n_i$ corners (so $\sum
  n_i=g$) then we can choose an $n_i$-fold branched cover $v_i\co
  D_i\to [0,1]\times\RR$ with the specified boundary behavior. Choose
  a $g$-fold branched cover $w\co E\to S^2$. (This is where we use the
  assumption that $g>1$.) Splicing together
  $\amalg_i v_i$ and $w$ gives a $g$-fold branched cover $u_{\bD,4}\co
  S_4\to[0,1]\times\RR$ with the desired boundary behavior. Perturbing
  $u_{\bD,4}$ slightly, we can ensure that all branch points of
  $u_{\bD,4}$ are simple.

  Consider the map $u=u_{\Sigma,4}\times u_{\bD,4}\co S_4\to
  \Sigma\times[0,1]\times\RR$. This satisfies
  Condition~(\ref{item:4.1p:bdy}),~(\ref{item:4.1p:one-arc})~(\ref{item:4.1p:components})
  and~(\ref{item:4.1p:acute}) of Lemma 4.1$'$, and the projection to
  $[0,1]\times\RR$ is a branched covering. Isotoping $u_{\bD,4}$ we
  can ensure Condition~(\ref{item:4.1p:holo}) as well as making
  $u_{\bD,4}$ holomorphic near its branch points.
  Finally, deforming $u$ slightly we can ensure that it has only
  double point singularities (Condition~(\ref{item:4.1p:double-pts})).
\end{proof}

\begin{proof}[Proof of Proposition $4.8'$]
  Using the fact that $\mu$ is unchanged by stabilization
  (see~\cite[Remark 10.5]{OS04:HolomorphicDisks}),
  as is $e+n_\vx+n_\vy$ (obvious), we may assume that $g>1$. 
  Since 
  \[
  \mu([\Sigma])=2=e([\Sigma])+n_\vx([\Sigma])+n_\vy([\Sigma]),
  \]
  adding or subtracting copies of $[\Sigma]$ changes both sides of
  Formula~(\ref{eq:mu-equals}) in the same way. So, it suffices to
  prove Formula~(\ref{eq:mu-equals}) after
  replacing $A$ by
  $A+(n+1)[\Sigma]$ where $A+n[\Sigma]$ is positive. Let $u$ be the
  map given by Lemma 4.8$'$ in the homology class
  $A+(n+1)[\Sigma]$. Via the tautological correspondence (see, for
  instance,~\cite[Section 13]{Lipshitz06:CylindricalHF}), $u$
  corresponds to a map $\phi\co \bD^2\to \Sym^g(\Sigma)$ with the same
  domain as $u$. Rasmussen showed~\cite[Theorem
  9.1]{Rasmussen03:Knots} that
  \begin{equation}\label{eq:Rasmussen}
  \mu(A)=\Delta\cdot \phi + 2e(A),
  \end{equation}
  where $\Delta$ denotes the diagonal.

  In terms of $u$, the intersections of $\phi$ with the $\Delta$ arise
  in two ways:
  \begin{itemize}
  \item Branch points of $\pi_\bD\circ u$. Lemma 4.8$'$ guaranteed
    that these be order $2$ branch points, and that $u$ be holomorphic
    near each of them. It follows that each branch point corresponds
    to a positive, transverse intersection of $\phi$ and $\Delta$.
  \item Double points of $u$. Each positive (respectively negative)
    double point corresponds to a positive (respectively negative),
    degree $2$ tangency of $\phi$ and $\Delta$.
  \end{itemize}
  So, we have
  \begin{align}\label{eq:Delta-dot-phi}
  \Delta\cdot\phi &= \br(\pi_\bD\circ u)+2(d_+-d_-),\\
  \shortintertext{where $\br$ denotes the number of branch points.
  By the Riemann-Hurwitz formula,} 
  \chi(S)&=g\chi(\bD^2)-\br(\pi_\bD\circ u)\nonumber
  \end{align}
  so
  \begin{align*}
  \br(\pi_\bD\circ
  u)&=g-\chi(S)=g-(g-n_\vx(A)-n_\vy(A)+e(A)+2(d_+-d_-))\\
  &=n_\vx(A)+n_\vy(A)-e(A)-2(d_+-d_-).
  \end{align*}
  Combining this with Equations~\eqref{eq:Rasmussen}
  and~\eqref{eq:Delta-dot-phi} gives
  \begin{align*}
  \mu(A)&=n_\vx(A)+n_\vy(A)-e(A)-2(d_+-d_-)+2(d_+-d_-)+2e(A)\\
  &=e(A)+n_\vx(A)+n_\vy(A),
  \end{align*}
  as desired.
\end{proof}

\section{Smaller errata}\label{sec:smaller}
\begin{enumerate}
\item In Section 1, on page 959, Condition (J4) should read ``$J(\bdy/\bdy
  s)=\bdy/\bdy t$'', not ``$J(\bdy/\bdy t)=\bdy/\bdy s$'' as currently
  written. (Thanks to C.~Taubes for pointing out this mistake.)
\item In Section 3, pages 966--972, the paper considers the space 
  \[
  W^{p,d}_{k}((S,\bdy S),(W,C_\alpha\cup C_\beta)).
  \]
  (See, for instance, Definition 3.5, page 968.) Here, $W=\Sigma\times[0,1]\times\RR$. It is important that this space of $W^{p,d}_k$ maps be a Banach manifold; this is used, for instance, in the proof of Proposition 3.7. However, it is not clear that this space of maps is a Banach manifold, because $W$ has boundary.
  
  The easiest way to fix this is to replace $W$ by $\Sigma\times\RR\times\RR$ (but leave the boundary conditions $C_\alpha$ and $C_\beta$ unchanged). This larger space has the structure of a Banach manifold in an obvious way. Since the projection to $\RR\times\RR$ is holomorphic, the $0$-set of the $\overline{\partial}$-operator on the larger space of maps is, in fact, contained in the smaller space of maps, so this has no effect on the space of holomorphic curves under consideration. (Thanks to J.~Pardon for pointing out this mistake.)
\item Also in Section 3, in the definition of the universal moduli
  space $\mathcal{M}^\ell$, instead of considering the space of all
  almost complex structures on $S$ (which is infinite-dimensional),
  one should consider the moduli space of complex structures on $S$
  (which is finite-dimensional). (Otherwise, in the proof of
  Proposition 3.8, the fiber $\mathcal{M}$ of the projection
  $\mathcal{M}\to\mathcal{J}$ is the product of the desired moduli
  space with an infinite-dimensional space.) (Thanks to J.~Pardon for
  pointing out this mistake.)
\item In Section 14.2, on page 1071, the form $df\wedge
  dg+\star\thinspace df$ is only closed if the Morse function $f$ is
  harmonic. To guarantee the existence of a harmonic Morse function,
  puncture the $3$-manifold $Y$ at two points, and consider functions
  which approach $+\infty$ at one of the punctures and $-\infty$ at
  the other puncture. (Thanks to G.~Xu for pointing out this mistake.)
\end{enumerate}

\bibliographystyle{hamsplain}\bibliography{heegaardfloer}

\providecommand{\bysame}{\leavevmode\hbox to3em{\hrulefill}\thinspace}
\providecommand{\eprint}{\begingroup \urlstyle{rm}\Url}
\begin{thebibliography}{1}

\bibitem{Lipshitz06:CylindricalHF}
Robert Lipshitz, \emph{A cylindrical reformulation of {H}eegaard {F}loer
  homology}, Geom. Topol. \textbf{10} (2006), 955--1097,
  \eprint{arXiv:math.SG/0502404}.

\bibitem{LOT1}
Robert Lipshitz, Peter~S. Ozsv{\'a}th, and Dylan~P. Thurston, \emph{Bordered
  {H}eegaard {F}loer homology: {I}nvariance and pairing}, 2008,
  \eprint{arXiv:0810.0687}.

\bibitem{OS04:HolomorphicDisks}
Peter~S. Ozsv{\'a}th and Zolt{\'a}n {\relax Sz}ab{\'o}, \emph{Holomorphic disks
  and topological invariants for closed three-manifolds}, Ann. of Math. (2)
  \textbf{159} (2004), no.~3, 1027--1158, \eprint{arXiv:math.SG/0101206}.

\bibitem{Rasmussen03:Knots}
Jacob Rasmussen, \emph{Floer homology and knot complements}, Ph.D. thesis,
  Harvard University, Cambridge, MA, 2003, \eprint{arXiv:math.GT/0306378}.

\bibitem{Sarkar06:IndexTriangles}
Sucharit Sarkar, \emph{{M}aslov index formulas for {W}hitney {$n$}-gons}, 2010,
  \eprint{arXiv:math/0609673 v3}.

\end{thebibliography}
\end{document}